\documentclass{article}
\usepackage{maa-monthly_mod}
\usepackage{amsmath,amssymb,amsthm, amscd, mathtools}
\usepackage{pifont,mathrsfs}
\usepackage{array,lastpage}
\usepackage{enumerate,xspace,pifont,shadow}
\usepackage{mathrsfs}
\usepackage[T1]{fontenc}
\usepackage{mathtools}
\usepackage{caption}
\usepackage{subcaption}

\final

\DeclareMathOperator{\vol}{vol} 
\DeclareMathOperator{\lcm}{lcm}
\def\vec#1{\mathchoice{\mbox{\boldmath$\displaystyle\bf#1$}}
{\mbox{\boldmath$\textstyle\bf#1$}}
{\mbox{\boldmath$\scriptstyle\bf#1$}}
{\mbox{\boldmath$\scriptscriptstyle\bf#1$}}}

\theoremstyle{theorem} 
\newtheorem{theorem}{Theorem}

\theoremstyle{definition}
\newtheorem{definition}{Definition}
\newtheorem*{remark}{Remark}
\newtheorem{example}{Example}
\newtheorem*{question}{Question}

\providecommand{\abs}[1]{\lvert#1\rvert}
\providecommand{\floor}[1]{\left\lfloor#1\right\rfloor} 
\providecommand{\Z}{\mathbb{Z}} \providecommand{\R}{\mathbb{R}}
 
 \providecommand{\x}{\mathbf{x}}
\providecommand{\y}{\mathbf{y}}

\begin{document}

\title{A Plethora of Polynomials:\\A Toolbox for Counting Problems}
\markright{A Plethora of Polynomials}
\author{Tristram Bogart and Kevin Woods}

\maketitle

\begin{abstract}
A wide variety of problems in combinatorics and discrete optimization depend on counting the set $S$ of integer points in a polytope, or in some more general object constructed via discrete geometry and first-order logic. We take a tour through numerous problems of this type. In particular, we consider families of such sets $S_t$ depending on one or more integer parameters $t$, and analyze the behavior of the function $f(t)=\abs{S_t}$. In the examples that we investigate, this function exhibits surprising polynomial-like behavior. We end with two broad theorems detailing settings where this polynomial-like behavior must hold. The plethora of examples illustrates the framework in which this behavior occurs and also gives an intuition for many of the proofs, helping us create a toolbox for counting problems like these.
\end{abstract}

\section{Introduction.}
\label{sec:intro}
It's no surprise that geometry is a fertile source of polynomial functions. For example, if we take a bounded $d$-dimensional object, $B$, and dilate by a factor of $t$, then $\vol(tB)=\vol(B)t^d$ is a degree $d$ polynomial in $t$. It's more surprising that we can add considerable complication to these sets, taking us into the worlds of \emph{discrete geometry} and \emph{first-order logic}, and we still see the persistent appearance of polynomial-like functions. The framework that we will describe is broad, and it includes many problems that, at first glance, seem unrelated to geometry and logic. We will examine a number of these problems as we build intuition and a toolbox. As a teaser that we will return to later, consider the \emph{Frobenius problem}:

\begin{definition}[Frobenius problem]
\label{def:frob}
Given positive integers $a_1,\ldots, a_n$, let $S$ be the semigroup generated by the $a_i$, that is, the set of integers that can be written as sums of the $a_i$'s:
\[S=\{m\in\Z_{\ge 0}:\ \text{there exist } \lambda_1,\ldots,\lambda_n\in\Z_{\ge 0}\text{ with } m=\lambda_1a_1+\cdots+\lambda_n a_n\}.\]
If the $a_i$ are relatively prime, then all sufficiently large integers are in $S$, and so we can define $F(a_1,\ldots,a_n)$ to be the largest integer \emph{not} contained in $S$ (the \emph{Frobenius number}) and $g(a_1,\ldots,a_n)$ to be the number of positive integers not contained in $S$ (the number of \emph{gaps}). (See \cite{alfonsin2005} for a broad overview.)
\end{definition}

\begin{example}
Chicken McNuggets used to be sold in packs of 6, 9, or 20. If you want to eat $m\in\Z_{\ge0}$ McNuggets, can you buy some packs and get exactly $m$ McNuggets? One can check that the only $m$ for which you \emph{cannot} are the 22 integers 1, 2, 3, 4, 5, 7, 8, 10, 11, 13, 14, 16, 17, 19, 22, 23, 25, 28, 31, 34, 37, and 43. That is, $F(6,9,20)=43$ and $g(6,9,20)=22$.
\end{example}

Finding formulas for $F(a_1,\ldots,a_n)$ and $g(a_1,\ldots,a_n)$ seems fruitless in general \cite{Ramirez96}, and so people have looked at many special cases. One consequence (among many) of the machinery that we explain here is that there are polynomial-like formulas whenever the $a_i$ are polynomial functions of a single parameter $t$. (This was first proved in \cite{RouneWoods}, when the $a_i$ are linear functions of $t$, and in \cite{Shen2015} for general polynomials.)

\begin{example} By the end of Section \ref{sec:nonlinearity}, we will have all of the tools that we need to calculate that 
\[F(t,t+1,t+3)=\begin{cases}
\frac{1}{3}t^2+t-1&\text{if $t\equiv 0\pmod 3$,}\\
\frac{1}{3}t^2+\frac{2}{3}t-2&\text{if $t\equiv 1\pmod 3$,}\\
\frac{1}{3}t^2+\frac{1}{3}t-1&\text{if $t\equiv 2\pmod 3$}
\end{cases}\]
and
\[g(t,t+1,t+3)=\begin{cases}
\frac{1}{6}t^2+\frac{1}{2}t&\text{if $t\equiv 0\pmod 3$,}\\
\frac{1}{6}t^2+\frac{1}{2}t-\frac{2}{3}&\text{if $t\equiv 1\pmod 3$,}\\
\frac{1}{6}t^2+\frac{1}{2}t-\frac{2}{3}&\text{if $t\equiv 2\pmod 3$.}
\end{cases}\]
\end{example}

In the next sections, we build up a library of examples --- starting simple and ending with ones like the above --- of polynomial-like functions that describe sets in discrete geometry and logic. Many of these examples illustrate important tools in the general study of these phenomena. If the reader is itching for a spoiler at any time, Section \ref{sec:thms} gives a precise statement of the general theorems encompassing all of these examples.

\section{Polyhedra.}
We start with a classic.
\begin{example}
\label{ex:tri1}
The $t$th triangular number, $1+2+\cdots+t$, is $t(t+1)/2$.
\end{example}
We've seen this fact so many times that it may have lost its wonder, but it truly is wonderful! An excited student, seeing this for the first time, would want to generalize. What other wonderful results are there like this?

\begin{figure}
\begin{center}
\includegraphics[width=\textwidth]{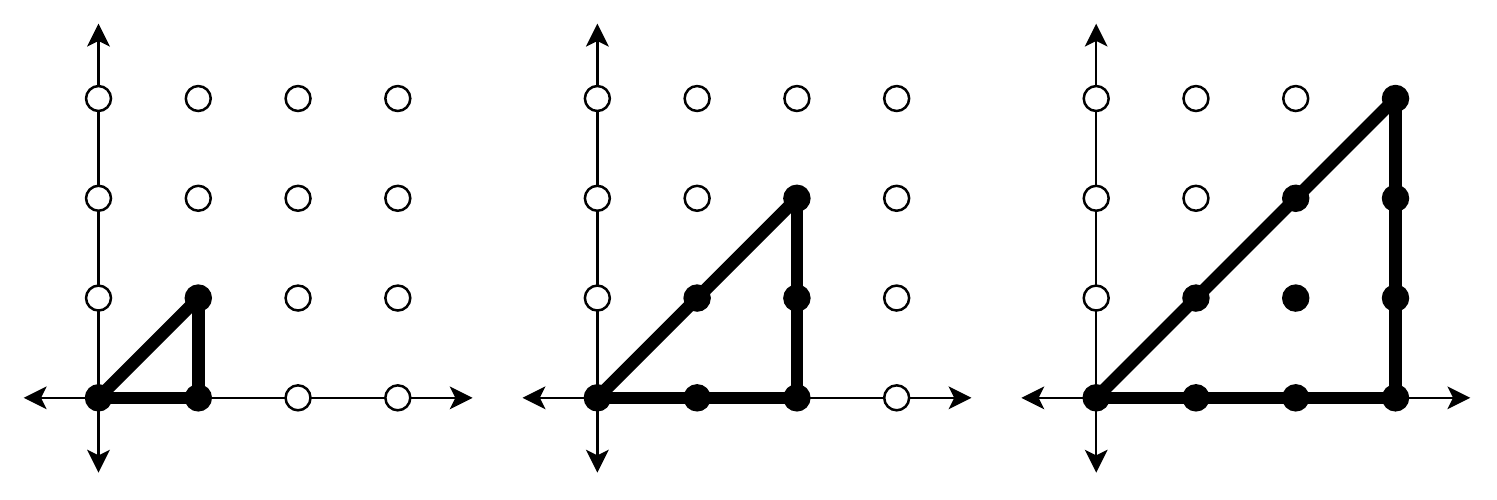}
\caption{$T\cap\Z^2$, $2T\cap\Z^2$, and $3T\cap\Z^2$ from Example \ref{ex:tri}, yielding the 2nd, 3rd, and 4th triangular numbers, respectively.}
\label{fig1}
\end{center}
\end{figure}

\begin{example}
\label{ex:tri}
Let $T\subseteq\R^2$ be the triangle with vertices $(0,0)$, $(1,0)$, and $(1,1)$, and let $tT$ (for $t\in\Z_{>0}$) be the dilation of $T$ by a factor of $t$, as pictured in Figure \ref{fig1}. Then we see that the number of points with integer coordinates that lie in $tT$ is the ($t+1$)st triangular number, so this is just the same example in a geometric disguise.
\end{example}

\begin{definition}
Given a bounded set $B\subseteq\R^d$, define $f_B(t)=\abs{tB\cap\Z^d}$, for $t\in\Z_{>0}$.
\end{definition}

For our triangle, $T$, in Example \ref{ex:tri}, we have
\[f_T(t)=\frac{(t+1)(t+2)}{2}=\frac{1}{2}t^2+\frac{3}{2}t+1.\]
It's not surprising that $f_T(t)$ is approximately $t^2/2$; after all, $\vol(tT)=t^2\vol(T)=t^2/2$, and so for large $t$ there should be roughly $t^2/2$ integer points contained in $tT$, with any discrepancy related to how the integer points interact with the boundary of $tT$. It is a beautiful fact that this discrepancy is so well behaved that $f_T(t)$ is still a polynomial.

How can we generalize? Let's start by examining a different triangle with integral vertices.

\begin{figure}
\begin{center}
\includegraphics[width=\textwidth]{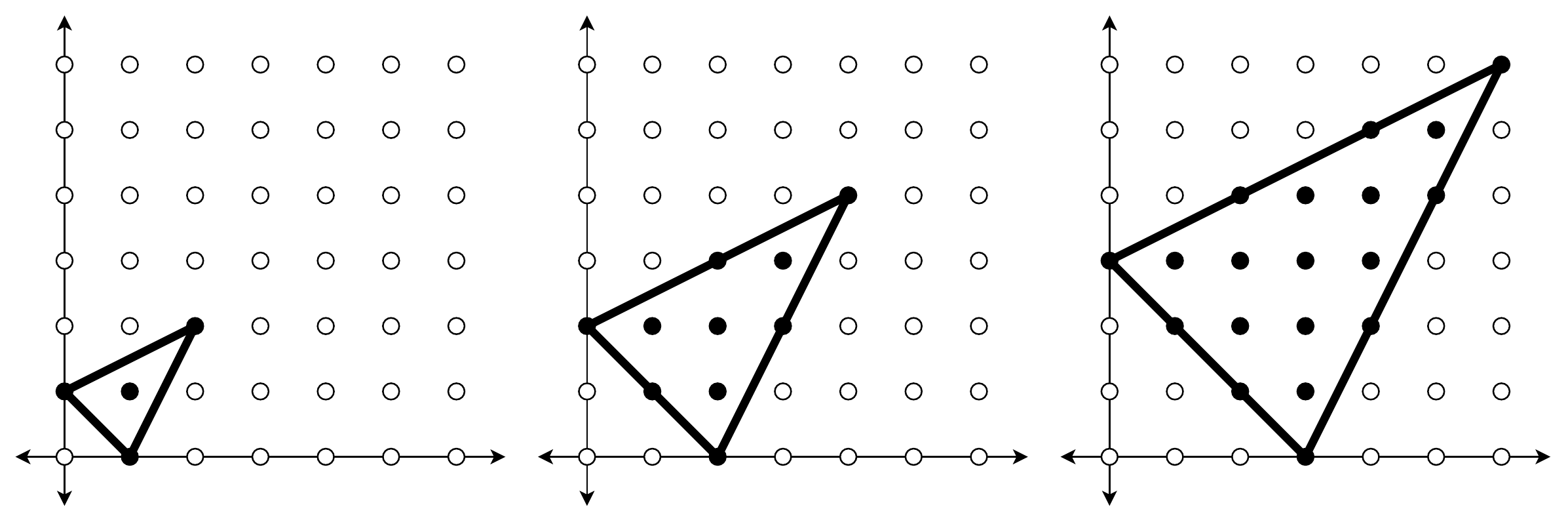}
\caption{$Q\cap\Z^2$, $2Q\cap\Z^2$, and $3Q\cap\Z^2$ from Example \ref{ex:triQ}.}
\label{fig2}
\end{center}
\end{figure}

\begin{example}
\label{ex:triQ}
Let $Q$ be the triangle with vertices $(1,0)$, $(0,1)$, and $(2,2)$; see Figure \ref{fig2}. Then
\[f_Q(t)=\frac{3}{2}t^2+\frac{3}{2}t+1.\]
\end{example}

It's fun to verify this formula by drawing dilates of $Q$, but we will establish it using Pick's theorem; in fact, this will prove that we get a nice polynomial formula for $f_P(t)$ when $P$ is \emph{any} polygon with integral vertices. Let $A_t$ be the area of $tP$, $b_t$ be the number of integer points on the boundary of $tP$, and $i_t$ be the number of integer points in the interior. Pick's theorem \cite{Pick} (see \cite[\S2]{BeckRobins} for a simple proof in English) tells us that
\[A_t=i_t+\frac{b_t}{2}-1.\]
Rearranging, we have that
\[f_P(t)=i_t+b_t=A_t+\frac{b_t}{2}+1.\]
But $A_t$ and $b_t$ are easy to calculate: $A_t=t^2A_1$ and $b_t=tb_1$ (the latter is easiest to see for the triangle $Q$ by dividing the boundary of $tQ$ into three half-open intervals, including one endpoint and excluding the other, each of which has $t$ vertices). Substituting, we have that
\[f_P(t)=A_1t^2+\frac{b_1}{2}t+1\]
is a polynomial. For the triangle $Q$, we use  $A_1=\frac{3}{2}$ and $b_1=3$, but we see that it works for any polygon with integral vertices. As before, it is not surprising that the leading term is $A_t=A_1t^2$; what is surprising is that the ``correction'' terms for looking at integer points are so simple.

How should we generalize next? There are two obvious ways: try dimensions higher than 2 or try polygons with nonintegral vertices. Don't worry: we'll try both. First let's see an example in a higher dimension.

\begin{figure}
\begin{center}
\includegraphics[width=\textwidth]{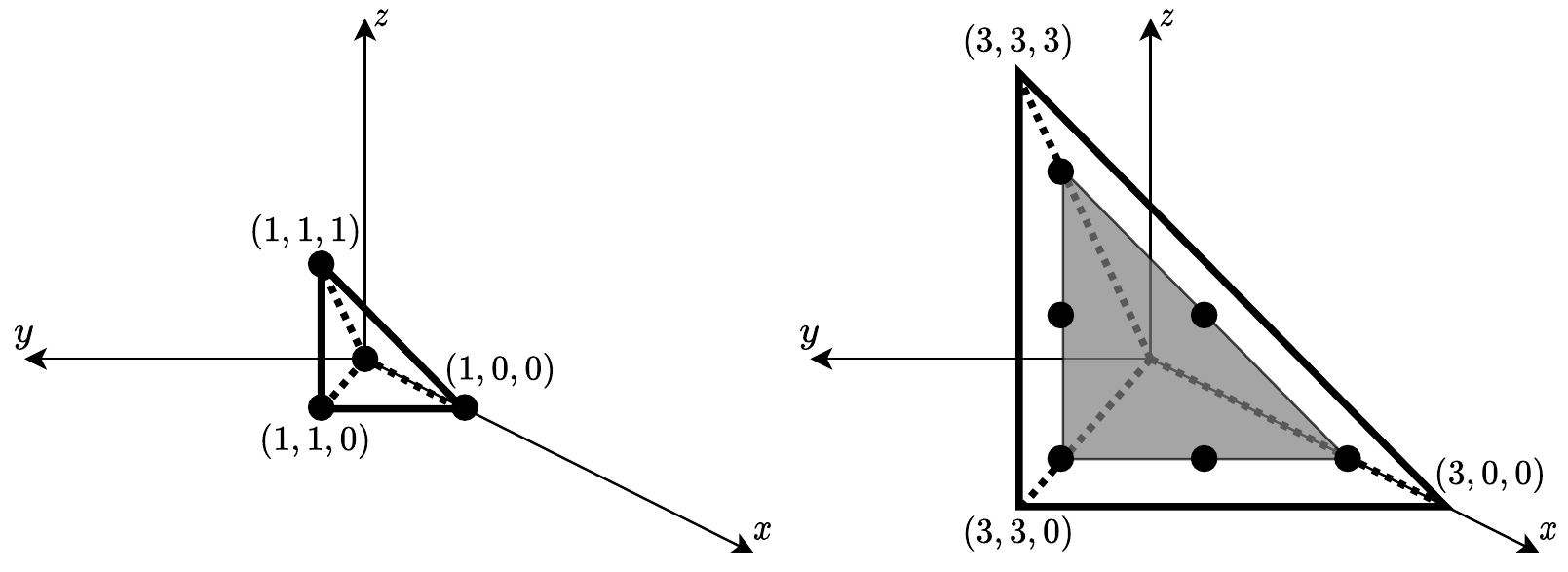}
\caption{On the left is $P$ from Example \ref{ex:tetra}, the tetrahedron formed by the convex hull of the three labeled vertices and the origin;  the $f_P(1)=\abs{P\cap\Z^2}=4$ integer points are shown. On the right is $3P$, and shaded is the cross-section of $3P$ at $x=2$; this cross-section has 6 integer points (the third triangular number), as shown.}
\label{fig3}
\end{center}
\end{figure}

\begin{example}
\label{ex:tetra}
Let $P\subseteq \R^3$ be the tetrahedron with vertices $(0,0,0)$, $(1,0,0)$, $(1,1,0)$, and $(1,1,1)$; see Figure \ref{fig3}. The cross-section of $tP$ at $x=s$ is a copy of $sT$, where $T$ is the triangle in Example \ref{ex:tri}. Thus the number of integer points in the cross-section is the $(s+1)$st triangular number, and so $f_P(t)$ is the sum of the first $t+1$ triangular numbers, that is, the $(t+1)$st \emph{tetrahedral number}. So we have
\[f_P(t)=\frac{(t+1)(t+2)(t+3)}{6}=\frac{1}{6}t^3+t^2+\frac{11}{6}t+1.\]
\end{example}

Again, our answer is a polynomial. Again the leading term is $\vol(tP)=\vol(P)t^3$. The first deep result in this area is that this works for \emph{any} polytope (bounded polyhedron) in \emph{any} dimension, as long as it has integral vertices. We'll state this precisely after a few examples with nonintegral vertices.

\begin{example}
Let $P\subseteq\R$ be the interval $[0,1/2]$. Then
\[f_P(t)=\floor{t/2}+1=\begin{cases} \frac{1}{2}t+1&\text{if $t$ even,}\\ \frac{1}{2}t+\frac{1}{2}&\text{if $t$ odd.}\end{cases}\]
\end{example}

This example makes it clear that, when the vertices are nonintegral, we can no longer expect to get exactly a polynomial. While $f_P(t)$ should be approximately $t/2$ (the length of $tP$), our formula will need to account for whether the right vertex of $tP$ is an integer ($t$ even) or not ($t$ odd). The most we can hope for is to get a counting function that is polynomial-\emph{like}, but includes some periodicity. Let's combine the polynomial behavior and periodic behavior into a definition.

\begin{definition}[Quasi-polynomial]
A function $f:\Z\rightarrow\Z$ is a \emph{quasi-polynomial} (which we will denote by QP) if there exist a period $s$ and polynomial functions $f_0(t), \ldots, f_{s-1}(t)$ such that $f(t)=f_i(t)$ for $t\equiv i\pmod s$.
\end{definition}

\begin{figure}
\begin{center}
\includegraphics[width=.8\textwidth]{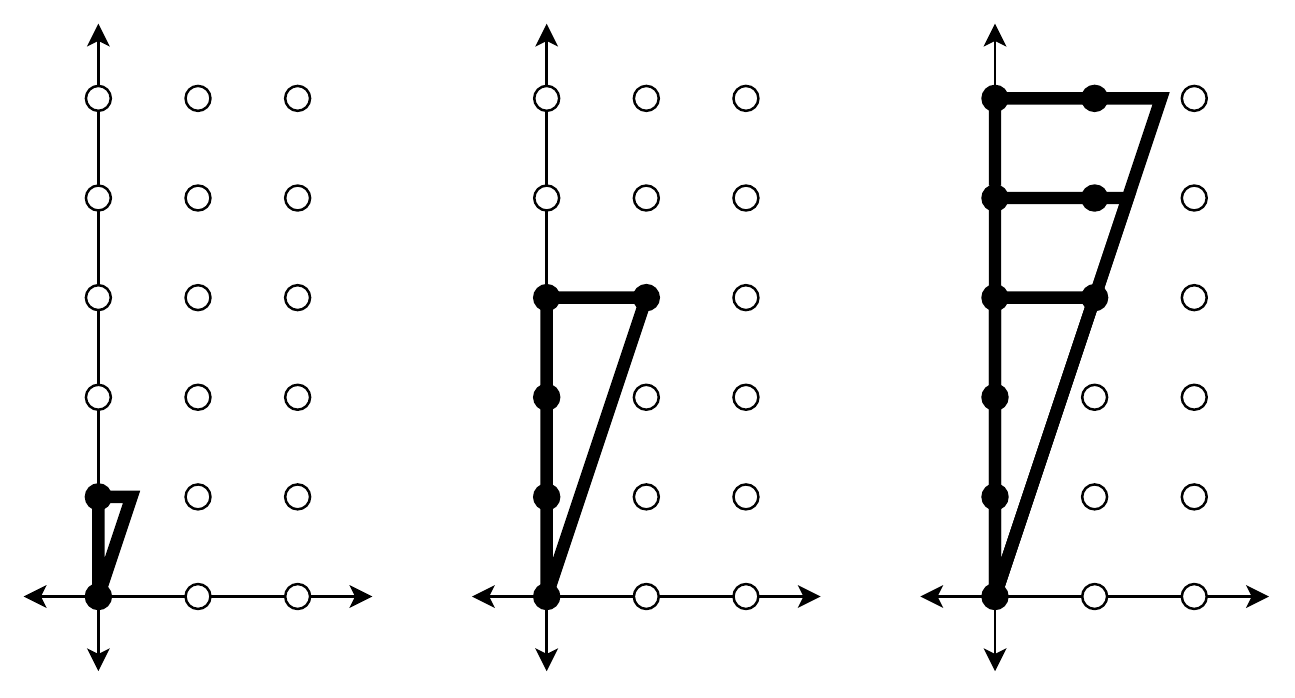}
\caption{$T\cap\Z^2$, $3T\cap\Z^2$, and an overlay of  $3T$, $4T$, $5T$, from Example \ref{ex:qp}. }
\label{fig4}
\end{center}
\end{figure}

\begin{example}
\label{ex:qp}
 Let $T$ be the triangle with vertices $(0,0)$, $(0,1)$, and $\left(\frac{1}{3},1\right)$; see Figure \ref{fig4}. How can we calculate $f(T)$? Observing that $3T$ is a polygon with integral vertices, area $3/2$, and containing five integer points on the boundary, we can use Pick's theorem to compute $f_T(t)$ when $t=3u$ is a multiple of three:
 \[f_T(t)=f_T(3u)=f_{3T}(u)=\frac{3}{2} u^2 + \frac{5}{2} u + 1= \frac{1}{6} t^2 + \frac{5}{6} t + 1.\]
On the other hand, if $t=3u+1$, then we can note that $(3u+1)T$ has exactly $u+1$ more integer points than $3uT$ has (they are $(0,3u+1),\ldots,(u,3u+1)$, to be precise; see the relationship between $3T$ and $4T$ in Figure \ref{fig4}). Therefore,
\[f_T(t)=f_T(3u+1)=f_{T}(3u)+(u+1)=\frac{1}{6} t^2 + \frac{5}{6} t + 1,\]
(coincidentally) the same polynomial as in the first case. We invite the reader to try this out for $t=3u+2$; successfully simplifying will yield
\[ f_T(t) = \begin{cases}
\frac{1}{6}t^2+\frac{5}{6}t+1&\text{if $t\equiv 0\pmod 3$,}\\
\frac{1}{6}t^2+\frac{5}{6}t+1&\text{if $t\equiv 1\pmod 3$,}\\
\frac{1}{6}t^2+\frac{5}{6}t+\frac{2}{3}&\text{if $t\equiv 2\pmod 3$}
\end{cases}\]
which is a QP of period three.
\end{example}

These results generalize to any polytope in any dimension, as long as it has rational vertices. This is the content of the first big theorem in this direction.

\begin{theorem}[Ehrhart's Theorem] \cite{Ehrhart62, Ehrhart77}.  Let $P\subseteq\R^d$ be a polytope with rational vertices, and let $s$ be the smallest positive integer such that $sP$ has integral vertices. Then $f_P(t)=\abs{tP\cap\Z^d}$ (defined for $t\in\Z_{>0}$) is a quasi-polynomial with period $s$. In particular, if $P$ has integral vertices, then $s=1$ and $f_P(t)$ is a polynomial.
\end{theorem}

Proofs of Ehrhart's theorem have been explained nicely in, for example \cite[\S 18]{barvinok2008} or \cite[\S 3]{BeckRobins}, so we do not delve into \emph{why} this is working in dimensions larger than two; we'd rather spend our time giving an intuition for the tools that layer on top of Ehrhart theory to prove QP results for larger classes of sets.

A final example shows that  polyhedra can be hidden in seemingly nongeometric problems.

\begin{figure}
\begin{center}
\includegraphics[width=.8\textwidth]{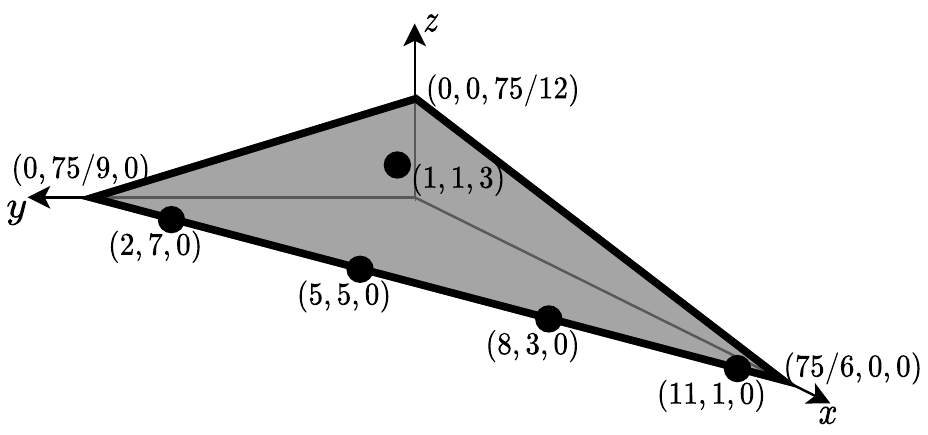}
\caption{The triangle $75P$, from Example \ref{ex:McNugget}. The three vertices are labeled, as are the $f_T(75)=5$ integer points in $75P$; for example, the point $(1,1,3)\in 75P$ shows that exactly 75 McNuggets can be purchased by ordering one box of 6, one box of 9, and three boxes of 20.}
\label{fig5}
\end{center}
\end{figure}

\begin{example}
\label{ex:McNugget}
Given that McNuggets come in boxes of 6, 9, and 20, let $f(t)$ be the number of different ways to order exactly $t$ McNuggets. For example, $f(18)=2$, since you can order three boxes of 6 or two boxes of 9. What is $f(t)$ as a function of $t$?

Let $T$ be the triangle, lying in three dimensions, with vertices $(1/6,0,0)$, $(0,1/9,0)$, and $(0,0,1/20)$. This triangle consists of the nonnegative real points lying in the hyperplane $6x+9y+20z=1$. So $tT$ is the triangle of nonnegative real points lying in the hyperplane $6x+9y+20z=t$, and an integer point in $tT$ is exactly a triplet of three nonnegative integers $(x,y,z)$ with $6x+9y+20z=t$; in other words, an integer point in $tT$ corresponds to a way to order exactly $t$ McNuggets. See Figure \ref{fig5} for a depiction of $75P$ and its 5 integer points. Ehrhart's theorem tells us that $f(t)=f_T(t)$ is a QP of period $\lcm(6,9,20)=180$ (since the vertices of $180T$ are integral). Calculating it is best done on a computer; for example, if $t\equiv 18\pmod{180}$, then $f(t)=\frac{1}{2160}t^2+\frac{7}{180}t+\frac{23}{20}$, and if $t\equiv 75\pmod{180}$, then $f(t)=\frac{1}{2160}t^2+\frac{11}{360}t+\frac{5}{48}$.
\end{example}

\section{Multiple parameters.}
Where to next? One direction to take is to think of a polytope, $P\subseteq\R^d$, as the solution set to a \emph{conjunction} of \emph{linear inequalities}. The dilate $tP$ can also be written as a conjunction of this sort, and it will have a very nice form: if $a_1x_1+\cdots +a_dx_d\le b$ is one of the defining inequalities of $P$, then $a_1x_1+\cdots +a_dx_d\le tb$ will be a defining inequality for $tP$. The facet of $tP$ that this inequality defines has a normal vector, $(a_1,\ldots, a_d)$, that is independent of $t$. However, the right-hand side of the inequality, $tb$, is changing, so the facets are shifting (but not twisting) as $t$ changes. The particular form of the right-hand side, $tb$, means that the facets shift in ``lockstep'' with each other, as linear functions of $t$. This suggests another way to generalize: allow the facets to shift separately from one another (but still with constant normal vectors), perhaps --- if we may be so bold --- by allowing additional parameters.

\begin{example}
\label{ex:trap}
Given $s,t\in\Z_{>0}$, define
\[P_{s,t}=\{(x,y)\in \R^2:\ y\ge 0,\ y\le s,\ y\le x,\ x+y\le t\}.\]
What is $f(s,t)=\abs{P_{s,t}\cap\Z^2}$?

\begin{figure}
\begin{center}
\includegraphics[width=\textwidth]{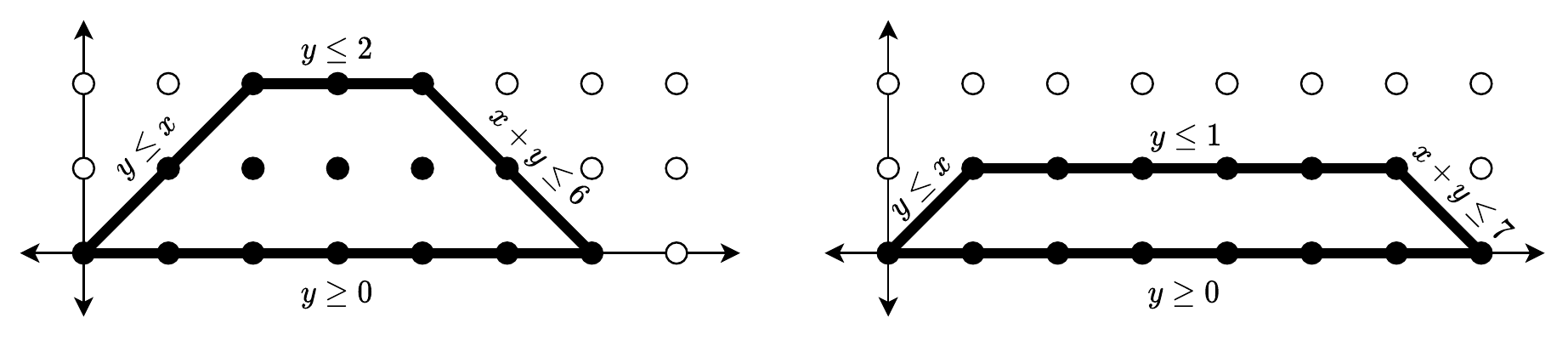}
\caption{$P_{2,6}\cap\Z^2$ and $P_{1,7}\cap\Z^2$  from Example \ref{ex:trap}.}
\label{fig6}
\end{center}
\end{figure}

Examining Figure \ref{fig6}, it looks like $P_{s,t}$ is a trapezoid, for fixed $s$ and $t$; as $s$ and $t$ change, two of the sides of this trapezoid shift. Counting the number of integer points isn't too bad: this trapezoid can be divided into two triangles and a rectangle, yielding total number of integer points
\[2\left(\frac{s(s+1)}{2}\right)+(t-2s+1)(s+1)=ts-s^2+t+1.\]
A polynomial! As before, if the vertices of $P_{s,t}$ were not integral, our answer would also depend periodically on $s$ and $t$. Perhaps you're shouting that we've missed something: if $s$ gets large but $t$ stays relatively small (to be precise, if $2s\ge t$), then the top edge of the trapezoid disappears, and we are left with a triangle with vertices $(0,0)$, $(t,0)$, and $(t/2,t/2)$. Therefore our final answer must be defined \emph{piecewise}, depending on whether $P_{s,t}$ is a trapezoid or a triangle. Putting it all together,
\[f(s,t)=\begin{cases} ts-s^2+t+1&\text{if $2s< t$,}\\
\frac{1}{4}t^2+t+1&\text{if $2s\ge t$, $t$ even,}\\
\frac{1}{4}t^2+t+\frac{3}{4}&\text{if $2s\ge t$, $t$ odd.}
\end{cases}\]
\end{example}

Functions like this are the final deviation from a nice, simple polynomial that we will need, made precise in the next two definitions.

\begin{definition}[Multivariate Quasi-polynomials]
\label{def:qp}
A function $f:\Z^k\rightarrow\Z$ is a (multivariate) \emph{quasi-polynomial} (which we will still denote by QP) if there exist a $k$-dimensional lattice $\Lambda \subseteq \Z^k$ and polynomial functions $f_C(t_1, \dots, t_k)$, one for each coset $C$ of $\Z^k$ modulo $\Lambda$, such that $f(t_1, \dots t_k) = f_C(t_1, \dots, t_k)$ if $(t_1, \dots, t_k) \in C$.
\end{definition}

\begin{remark}
In the $2s\ge t$ cases of the function $f(s,t)$ above, the lattice $\Lambda$ is $\Z\times 2\Z$ with cosets $(0,0)+\Lambda$ ($s$ any integer, $t$ even) and $(0,1)+\Lambda$ ($s$ any integer, $t$ odd). In general, $\Lambda$ might be a more complicated sublattice of $\Z^k$, but one could always define the quasi-polynomials on cosets of some rectangular lattice $m_1\Z \times \dots \times m_k\Z \subseteq \Lambda$ so that $f(t_1,\dots,t_k)$ is defined by polynomials on residue classes of $t_i$ mod $m_i$ for each $i = 1, \dots, k$.
\end{remark}

\begin{definition}[Piecewise Quasi-polynomials]
\label{def:pqp}
A function $f:\Z^k\rightarrow \Z$ is a \emph{piecewise quasi-polynomial} (PQP) if there is a partition of $\R^k$ into finitely many pieces $S_1, \dots, S_q$, each defined by a conjunction of linear inequalities, such that for each $i$, $f$ agrees with a multivariate quasi-polynomial on $S_i$.
\end{definition}



\begin{remark} Depending on the application, we often restrict the domain of our PQP to be, for example, $\Z^k_{\ge 0}$.
\end{remark}

\begin{remark} Geometrically, each piece $S_i$ is either a bounded polytope or an unbounded polyhedron. If $k=1$, the partition of $\R_{\ge 0}$ in the definition of a PQP $f:\Z_{\ge 0}\rightarrow \Z$ must contain exactly one unbounded polyhedron, and so $f$ must agree with a quasi-polynomial for all sufficiently large $t$; therefore, we call $f$ an \emph{eventual quasi-polynomial} (EQP).
\end{remark}

\begin{example}
McDonald's burger options are a little confusing. You can buy a Double Cheeseburger, which has two patties and two slices of cheese. Or you can buy a McDouble, which is identical except it has only one slice of cheese. The difference in price is completely unpredictable. Of course, you could buy a Triple Cheeseburger, which has three patties and (obviously?) two slices of cheese. Feeling particularly peckish, you want to consume 14 patties and 9 slices of cheese. Can you order exactly this (ignoring the buns), using these three types of sandwiches? In how many ways? We leave it to the reader to find the 3 different ways this can be done and to compute the number of ways to order $s$ patties and $t$ slices of cheese, in general. The answer is a PQP! For example, if $s\le 2t$, $2s\ge 3t$, and $s$ is even, then there are $(2t-s+2)/2$ ways.

Again, this problem can be written as a conjunction of linear inequalities: if $x,y,z\in\Z$ are the number of Double Cheeseburgers, McDoubles, and Triple Cheeseburgers, respectively, then we must have $x,y,z\ge 0$, $2x+2y+3z=s$, and $2x+1y+2z=t$ (noting that an equality $=$ is just the conjunction of two inequalities $\ge$ and $\le$).
\end{example}

This example is a special case of the \emph{vector partition function}, a multi-dimensional generalization of a Frobenius-style problem: for fixed vectors $\vec v_1,\ldots,\vec v_n\in\Z^k_{\ge 0}$, given $\vec t\in\Z^k_{\ge 0}$, define $f(\vec t)$ to be the number of nonnegative integer vectors $(\lambda_1,\ldots,\lambda_n)$ such that $\sum_i\lambda_i\vec v_i=\vec t$ (we use bold notation to signify a vector). Then $f(\vec t)$ will be a PQP \cite{sturmfels1995}. In the example, we have $\vec v_1=(2,2)$, $\vec v_2=(2,1)$, and $\vec v_3=(3,2)$.

\subsection{Where we stand.}
The sets we are defining have two types of variables, the \emph{parameter} variables $\vec t\in\Z^k$ and the \emph{counted} variables $\vec x\in\Z^n$. We define a parametric family of sets $S_{\vec t}$, and examine $f(\vec t)=\abs{S_{\vec t}}$. In the examples so far, $S_{\vec t}$ is the set of $\vec x\in\Z^n$ such that the variables $\vec t, \vec x$ satisfy a conjunction of linear inequalities (interestingly, $\vec t$ and $\vec x$ are indistinguishable in the linear inequalities; they are only distinguished when we decide to fix $\vec t$ and look at the set, $S_{\vec t}$, of $\vec x$ that these inequalities define). In this setup, we have seen examples where $f(\vec t)=\abs{S_{\vec t}}$ is a PQP. This will always be true! See \cite[\S 18]{barvinok2008}, for example, where this is proved alongside Ehrhart's theorem.



\section{Boolean Combinations.}
Where to next? Once we start thinking of a polyhedron as a \emph{conjunction} of linear inequalities, it seems natural to allow arbitrary boolean combinations. That is, we allow our linear inequalities to be joined together with any combination of logical operators: $\wedge$'s (and's), $\vee$'s (or's), and $\neg$'s (not's). This gives us some great examples of QP counting functions, many of which are not obviously geometric in nature.

\begin{figure}
\begin{subfigure}{.5\textwidth}
  \centering
  \includegraphics[width=.8\linewidth]{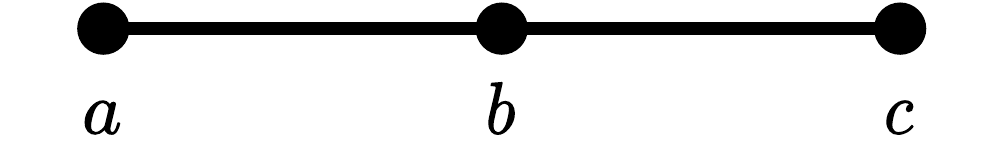}  
  \caption{The graph, $G$, in Example \ref{ex:chrom}.}
  \label{fig7}
\end{subfigure}
\begin{subfigure}{.5\textwidth}
  \centering
  \includegraphics[width=.6\linewidth]{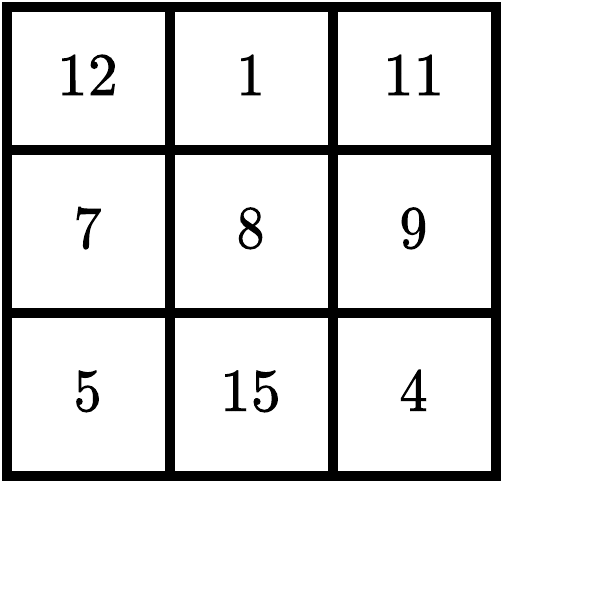}  
  \caption{A 24-magic square; see Example \ref{ex:mag}.}
  \label{fig7b}
\end{subfigure}
\caption{\ }
\end{figure}

\begin{example}
\label{ex:chrom}
Let $G$ be the graph with vertices $a,b,c$ and two edges $\{a,b\}$ and $\{b,c\}$; see Figure \ref{fig7}. The \emph{chromatic polynomial}, $\chi_G(t)$, is the number of ways to color the vertices of $G$ with $t$ possible colors, so that no adjacent vertices have the same color. For this particular graph $G$, calculating that $\chi_G(t)=t(t-1)^2$ is a standard exercise. (Hint: how many ways can you color $a$? Now $b$? Now $c$?) But it also fits our current setup of boolean combinations of linear inequalities: if $(x_a,x_b,x_c)$ are the colors of $a$, $b$, and $c$, respectively, with colors corresponding to the integers $1,2,\ldots,t$, then a coloring of $G$ satisfies
\[(1\le x_a,x_b,x_c\le t)\ \wedge\ \neg(x_a=x_b) \ \wedge\ \neg(x_b=x_c).\]
\end{example}

Chromatic polynomials, along with many other interesting examples, were first conceived in this way in \cite{beck2006inside} as the counting functions of \emph{inside-out polytopes}.

%

\begin{example}
\label{ex:mag}
Define a $t$-magic square to be a way to fill a $3\times 3$ grid with nine distinct positive integers such that the sum of every row, every column, and the two diagonals are all exactly $t$; see Figure \ref{fig7b} for a $24$-magic square. The number of $t$-magic squares depends on $t\pmod{18}$ (see \cite{beck2010six} for this and many similar calculations); if $t$ is not a multiple of 3, then there are no $t$-magic squares (can you prove it?).  On the other hand, for example, if $t\equiv 6\pmod{18}$, then there are $\frac{2}{9}(t-6)(t-10)$ such $t$-magic squares. Again a quasi-polynomial!

If $x_{ij}$ is the number in the $i$th row and $j$th column, then we see that the set of $t$-magic squares can be defined as a boolean combination of linear inequalities. For example, we require $x_{11}\ne x_{12}$ and so on (all entries distinct), and $x_{11}+x_{12}+x_{13}=t$ and so on (all rows, columns, diagonals sum to $t$).
\end{example}

\begin{example}
How many ways are there to place three queens on an $t\times t$ board such that no two queens are attacking each other?  This and many other formulas are calculated in \cite{chaiken2019q}: there are
\[\frac{t^6}{6}-\frac{5t^5}{3}+\frac{79t^4}{12}-\frac{25t^3}{2}+11t^2-\frac{43t}{12}+\frac{1}{8}+(-1)^t\left(\frac{t}{4}-\frac{1}{8}\right)\]
ways, a QP of period 2. This problem falls into our setting! Let's momentarily label the queens 1, 2, 3 and count labeled ways to place them (and then divide by 6 since all permutations of the three queens are equivalent). Suppose we place Queen $i$ at position $(x_i,y_i)\in\Z^2$ where $1\le x_i,y_i\le t$. The nonattacking conditions are given by boolean combinations of linear inequalities; for example, $x_1\ne x_2$ says that the first two queens can't be in the same row, and $x_1-y_1\ne x_2-y_2$ says they cannot be along one of the same diagonals.
\end{example}

\subsection{Where we stand.} We can define our sets $S_{\vec t}$ using boolean combinations of linear inequalities in the parameters $\vec t$ and counted variables $\vec x$. If $f(\vec t)$ is the number of $\vec x\in S_{\vec t}$, then $f(\vec t)$ will be PQP. Allowing boolean combinations opens up a wide variety of possible applications of these methods. Mathematically, however, there is not much new: boolean combinations may be converted into \emph{disjunctive normal form}, that is, written as a disjunction of conjunctions of linear inequalities, that is, a union of polyhedra. With a little care on the overlaps, we can rewrite such a set as a \emph{disjoint} union of polyhedra, and therefore apply the old methods.

\section{Quantifiers.}
\label{sec:quant}
Where to next? Having allowed boolean combinations, we've set our list of allowed operations firmly in a logical realm. Adding quantifiers ($\exists x\in\Z$ and $\forall x\in\Z$) gives us a full first-order logic, such an important one that it has a name: \emph{Presburger arithmetic} \cite{Presburger29} (see \cite{Presburger91} for a translation). The most obvious thing that allowing quantifiers buys us is that we may now define sets with regular gaps in them, such as the set of odd numbers:

\begin{example}
Given $t\in\Z_{\ge 0}$ define
\[S_t=\{x\in\Z:\ \exists y\in\Z,\ x=2y+1\ \wedge\ 0\le x\le t\},\]
the set of odd integers between 0 and $t$. Happily, we still get quasi-polynomial behavior: $\abs{S_t}=\floor{(t+1)/2}$.
\end{example}

Allowing quantifiers definitely gives us a robust collection of sets we can now describe, and the good news is that we will still get quasi-polynomial behavior. The key is a process called \emph{quantifier elimination}, explained nicely in \cite{Cooper72}.
For simplicity, we illustrate this method with an example that doesn't have any parameters. 

\begin{figure}
\begin{center}
\includegraphics[width=.6\textwidth]{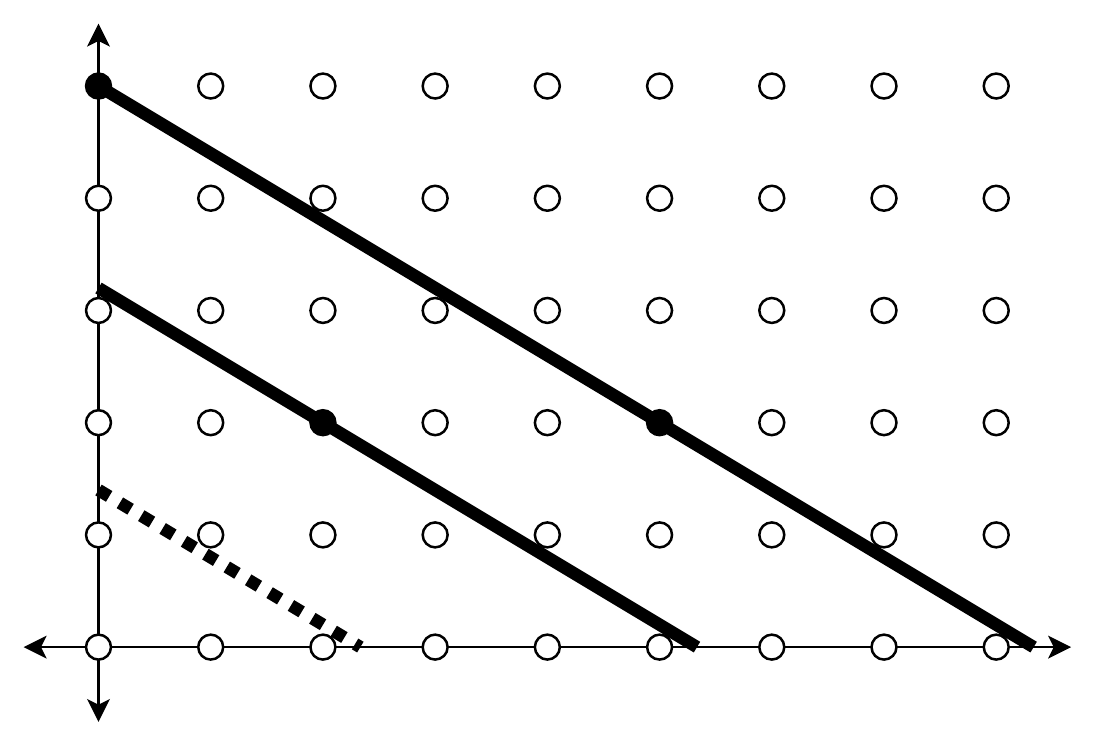}
\caption{Diagonal line segments are $3x+5y=n$, for $n=7,16,25$, constrained by $x,y\ge 0$. We see that $7\notin S$, because there are no integer points on the line segment. On the other hand $16\in S$ and $25\in S$; in fact, the two integer points for $n=25$ correspond to two ways to write 25 using 3's and 5's: $25=0\cdot 3+5\cdot 5=5\cdot 3+2\cdot 5$. Furthermore, we can see the second quantifier elimination step used in Example \ref{ex:qe}: if there are some $x,y\in\Z_{\ge 0}$ with $3x+5y=n$, then we may take $y$ to be 0, 1, or 2; indeed, if $(x,y)$ is on the line segment with $y\ge 3$, then $(x+5,y-3)$ is also on the line segment, e.g., $(0,5)$ and $(5,2)$ are on the line segment corresponding to $n=25$. }
\label{fig8}
\end{center}
\end{figure}

\begin{example}
\label{ex:qe}
Suppose we want to understand the semigroup generated by 3 and 5,
\begin{align*}S&=\{n\in\Z:\ \exists x,y\in\Z, \ (x\ge 0)\ \wedge\ \ (y\ge 0)\ \wedge\ \ (3x+5y=n)\}\\
&=\{0,3,5,6,8,9,10,\ldots\};\end{align*}
see Figure \ref{fig8}.

The quantifiers are the new complication here; we'd like to get rid of them, that is, describe $S$ \emph{without} needing to use quantifiers at all. Quantifier elimination is a tool for doing this. We first get rid of the existential quantifier on $x$. This is pretty straightforward: if there is any integer $x$ that works, it must be $x=(n-5y)/3$, to make the equation $3x+5y=n$ true. We need to be a little careful here: $x$ must be an integer, so we must have $3\big| (n-5y)$. In that case, we can simply substitute $x=(n-5y)/3$ into the constraints: $x\ge 0$ becomes $(n-5y)/3\ge 0$, $y\ge 0$ stays the same, and $3x+5y=n$ becomes the tautology $n=n$. After simplifying slightly, we have found a new way to express the semigroup:
\[S=\Big\{n\in\Z:\ \exists y\in\Z,\ 3\big| (n-5y)\ \wedge\ (5y\le n)\ \wedge \ (y\ge 0)\Big\}.\]
At the minor cost of introducing the \emph{divisibility predicate} $3\big| (n-5y)$, the quantified variable $x$ has been eliminated!

Let's try to eliminate $y$ next. This is not quite as easy, since there is no equality that would say exactly what $y$ must be. We instead ask ourselves, ``If there exists \emph{some} $y$ that makes this true, what is the smallest $y$ that works?'' We claim that the smallest $y$ must be 0, 1, or 2: indeed, exactly one of these will make $3\big| (n-5y)$ true, these are the smallest $y$'s that might make $y\ge 0$ true, and to make $5y\le n$ we actually \emph{want} to choose $y$ as small as possible (this is, if $y=i$ satisfies $5y\le n$, then so does $i-3$, $i-6$, etc., and so either $y=0$, $y=1$, or $y=2$ will satisfy it). In other words, we may replace this expression for $S$ with one that has three cases (separated with or's), substituting $y=0$, $y=1$, and $y=2$, respectively:
\begin{align*}\Big(3\big|& n\ \wedge \ (0\le n)\ \wedge\ (0\ge 0)\Big)\ \vee\ \Big(3\big| (n-5)\ \wedge \ (5\le n)\ \wedge\ (1\ge 0)\Big)\\
&\vee\ \Big(3\big| (n-10)\ \wedge \ (10\le n)\ \wedge\ (2\ge 0)\Big).\end{align*}

The quantifiers are gone, leaving only the variable $n$. Wow! We can now apply all of the machinery that we had already built up for quantifier-free formulas. In this particular example, we have divided into three disjoint cases based on whether $n$ is 0, 2, or 1 modulo $3$, respectively:
\[S=\{0,3,6,\ldots\}\sqcup\{5,8,11,\ldots\}\sqcup\{10,13,16,\ldots\},\]
and the set of nonnegative integers \emph{not} in $S$ is
\[ \emptyset \sqcup\{2\} \sqcup\{1,4,7\}.\]
Thus $F(3,5)=7$ and $g(3,5)=4$ (see Definition \ref{def:frob}). In fact, the pattern of quantifier elimination is fairly clear if we replace 3 and 5 by any relatively prime $a$ and $b$: it will give us the decomposition of $S$ based on its so-called \emph{Ap\'ery sets} \cite{Apery} (and see \cite{Ramirez96} for much more):
\[S=\bigsqcup_{i=0}^{a-1}\left( bi+a\Z_{\ge 0}\right)\qquad \text{and} \qquad S^c=\bigsqcup_{i=0}^{a-1} \{bi-aj:\ 1\le j\le\floor{bi/a}\}.\]
The maximum element of $S^c$ occurs when $i=a-1$ and $j=1$, and so  $F(a,b)=b\cdot (a-1)-a\cdot 1=ab-a-b$. Less obviously, we also get $g(a,b)=(ab-a-b+1)/2$. (Hint: pair the $i$ and $a-i$ terms in this union together.) These formulas were first proved in \cite{Sylvester84}.
\end{example}

\begin{remark} The number of gaps $g(a,b)$ is a counting function as in all of our other examples, but the Frobenius number $F(a,b)$ has a slightly different character: it \emph{specifies} the maximal element of the set of gaps. In general, we might want to do something like specify the maximal element of any of our families of one-dimensional sets $S_t$, or in higher dimensions specify an element that maximizes a given linear functional. The general results for specification problems are similar to those for counting problems, as we will see in Section \ref{sec:thms}.
\end{remark}

Now we give an example that does have a parameter (three parameters!), for which we get a PQP counting function. This example could be analyzed with the same quantifier elimination tools as above (we will leave the reader to have fun with it).

\begin{example}
A particular McDonald's location will have only a finite number of boxes of McNuggets available to buy. Suppose they have $r$ boxes of 6, $s$ boxes of 9, and $t$ boxes of 20 available. How many different numbers of McNuggets could you buy? The answer will be PQP, since the set of possibilities can be defined with quantifiers, boolean combinations, and linear inequalities:
\begin{align*}S_{r,s,t}=\{n&\in\Z:\ \exists x,y,z\in\Z,\\
&(0\le x\le r)\  \wedge\ (0\le y\le s)\ \wedge\ (0\le z\le t)\ \wedge\ (n=6x+9y+20z)\}.\end{align*}
Indeed, for $r$, $s$, and $t$ sufficiently large, one can verify that there will be $6r+9s+20t-43$ possibilities.
\end{example}

For an algebraic perspective that interprets functions like this as Hilbert polynomials, see \cite{khovanskii1992newton}.

\subsection{Where we stand.} We can define our sets $S_{\vec t}$ using quantifiers and boolean combinations of linear inequalities; there will now be three types of variables: the parameters $\vec t$, the counted variables $\vec x$, and the \emph{bound} variables $\vec y$ (those associated with a quantifier).  If $f(\vec t)$ is the number of $\vec x\in S_{\vec t}$, then $f(\vec t)$ will be PQP. The main tool is quantifier elimination, which reduces the problem to the previously understood quantifier-free case; see \cite{woods2015presburger}.

\section{Nonlinearity, with one parameter variable.}
\label{sec:nonlinearity}
Where to next? For our last twist (pun intended), we must restrict ourselves to only one parameter $t$. Recalling prior polyhedral examples, changing $t$ only caused a parallel shift of the facets of the polyhedra. Now we will allow the facets to ``twist'' by allowing the normal vectors to also depend on $t$.

\begin{figure}
\begin{center}
\includegraphics[width=.5\textwidth]{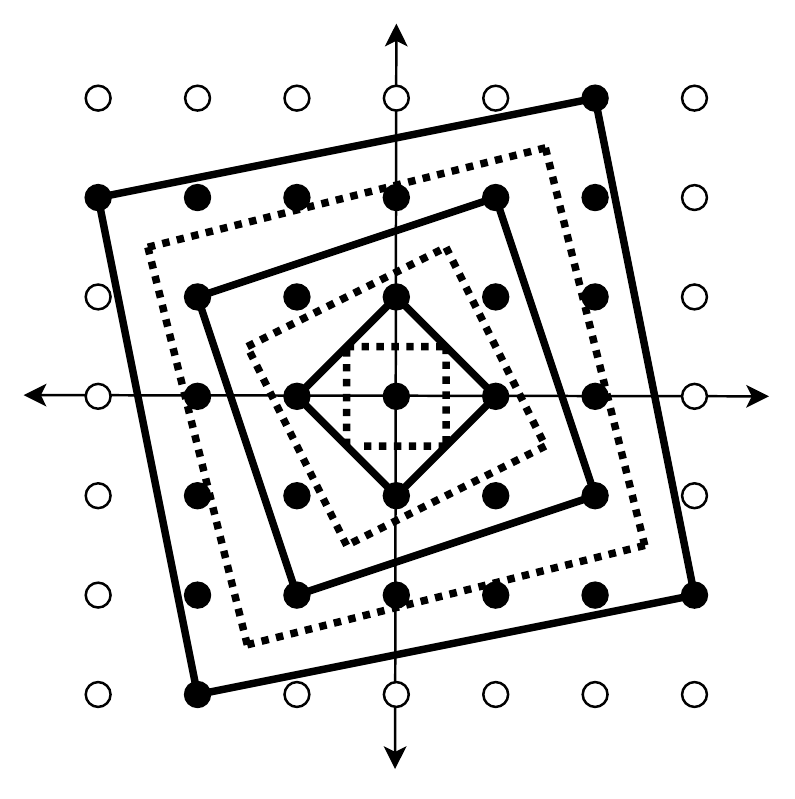}
\caption{The twisting square $P_t$ in Example \ref{ex:twist}, for $t=1,2,3,4,5,6$. Even $t$ are in solid lines and odd $t$ in dotted lines.}
\label{fig9}
\end{center}
\end{figure}

\begin{example} \cite{Woods_unreasonable}
\label{ex:twist}
Let $P_t$ be the ``twisting square'' defined by
\begin{align*}P_t=\Big\{(x,y)\in\mathbb R^2:\  &-(t^2-2t+2)\le2x+(2t-2)y\le t^2-2t+2,\\
&-(t^2-2t+2)\le (2-2t)x+2y\le t^2-2t+2\Big\};
\end{align*}
see Figure \ref{fig9}. Then $\abs{P_t\cap\Z^2}$ is given by the QP
\[\abs{P_t\cap\Z^2}=\begin{cases}t^2-2t+2 &\text{if $t$ odd,}\\t^2-2t+5 &\text{if $t$ even}. \end{cases}
\]
We leave the proof that this formula for $\abs{P_t\cap\Z^2}$ is correct as an exercise.  
\end{example}

Remarkably, lattice-point counting functions are EQP (eventual quasi-polynomials, that is, quasi-polynomial for sufficiently large $t$) for all such families of twisting polytopes \cite{CLS}. To be precise, this result applies to polytopes that can be written as conjunctions of inequalities of the form
 \[ a_1(t) x_1+\cdots +a_d(t) x_d\le b(t) \]
  where $a_1(t), \dots, a_d(t), b(t)$ are polynomials in $t$ with integer coefficients. That is, we no longer require our defining inequalities to be linear; we may now multiply by $t$ with impunity. Indeed, we can do this, not only with polytopes, but with any formula built from quantifiers, boolean operations, and inequalities of the above form. This was conjectured in \cite{Woods_unreasonable}, proved in \cite{BGW2017}, and will be stated precisely in the next section. We first proceed with some examples.
  
For these first two examples, note that two very natural functions, integer division and gcd, can be expressed in this new language: if $f(t)$ and $g(t)$ are polynomials with integer coefficients, then $\floor{f(t)/g(t)}$ is the largest $x\in\Z$ such that $g(t)x\le f(t)$, and $\gcd\big(f(t),g(t)\big)$ is the smallest $y\in\Z_{>0}$ such that
\[\exists a,b\in\Z,\ af(t)+bg(t)=y.\]

The next two examples illustrate how to see that $\floor{f(t)/g(t)}$ and $\gcd\big(f(t),g(t)\big)$ are indeed EQP functions.
\begin{example}

Let's evaluate $\floor{\frac{t^2}{2t+3}}$. By polynomial division, we obtain
  \[ \frac{t^2}{2t+3} = \frac{t}{2}-\frac{3}{4} + \frac{9/4}{2t+3}=\frac{t-2}{2}+\left(\frac{1}{4}+ \frac{9/4}{2t+3}\right). \]
 For sufficiently large $t$, the term in parentheses is strictly between 0 and $\frac{1}{2}$, and thus
  \[ \floor{\frac{t^2}{2t+3}}=\floor{\frac{t-2}{2}}=\begin{cases} \frac{t-2}{2} \text{ if $t \equiv 0 \pmod{2}$,} \\
      \frac{t-3}{2} \text{ if $t \equiv 1 \pmod{2}$}\end{cases} \]
  for all $t \geq 4$ (but not for $t=1$ or $t=3$).
\end{example}

%
%
%

\begin{example}
Let's calculate $\gcd\left(t^2+1,2t-1\right)$. A natural way to do this is to use the Euclidean algorithm, and the previous example shows that the quotients and remainders are all EQPs, so this algorithm must work to calculate the gcd as an EQP. In this example, we divide $2t-1$ into $t^2+1$ and immediately see that our answer must depend on $t$ mod 2. Let's look at even $t$, so $t=2u$ for some $u\in\Z$, $t^2+1=4u^2+1$, and $2t-1=4u-1$. The Euclidean algorithm yields
\begin{align*}4u^2+1&=u(4u-1)+(u+1)\\
4u-1&=4(u+1)-5,
\end{align*}
and so
\[\gcd(4u^2+1,4u-1)=\gcd(u+1,5)=\begin{cases}5 & \text{if $u\equiv 4\pmod{5}$,}\\ 1 & \text{otherwise.}\end{cases}\]
Note that $u\equiv 4\pmod 5$ corresponds to $t\equiv 2\cdot 4\equiv 3\pmod 5$. We leave the reader to try the case where $t$ is odd, which turns out to give the same answer: regardless of parity, we have
\[\gcd\left(t^2+1,2t-1\right)=\begin{cases} 5&\text{if $t\equiv 3\pmod 5$,}\\ 1&\text{otherwise}.\end{cases}\]
\end{example}

 Let's finally tackle our teaser problem from Section \ref{sec:intro}, to compute $F(t,t+1,t+3)$ and $g(t,t+1,t+3)$, the Frobenius number and number of gaps for the semigroup generated by $t$, $t+1$, and $t+3$, as defined in Definition \ref{def:frob}. These problems involve most of the twists we have introduced; in particular, they require both quantifiers and multiplication of the parameter variable $t$ by other variables. 

\begin{example}
Define the semigroup $S_t$ to be the set of integers $n$ such that
\[\exists x,y,z\in \Z_{\ge 0},\ tx+(t+1)y+(t+3)z=n.\]
We'll just do the first step for calculating $F(t,t+1,t+3)$ and $g(t,t+1,t+3)$, which requires a new technique, the ``base $t$'' method from \cite{CLS}. The remainder of the computation will be doable with techniques from previous sections, but it would be calculation intensive. 
Note that for all $n\ge t^2$, we have $n\in S_t$ (indeed $F(t,t+1)=t(t+1)-t-(t+1)=t^2-t-1$ already, so everything larger than $t^2-t-1$ must be in $S_t$). So we concentrate on $n<t^2$, and we write $n$ in ``base $t$,'' that is, $n=at+b$ with $0\le a,b<t$. Also note that $0\le x,y,z<t$ must be true in our defining equation, $tx+(t+1)y+(t+3)z=n$ (or else $n\ge t^2$). We rearrange our defining equation:
\[t(x+y+z)+(y+3z)=n=at+b.\]
Note that if $y+3z<t$, then the left-hand side must be a ``base $t$ expression'' (a polynomial in $t$ with coefficients between $0$ and $t-1$); this would mean that our expression simplifies to
\[(x+y+z=a)\ \wedge (y+3z=b).\]
Of course, $y+3z$ need not be less than $t$, but we do know it is less than $t+3t=4t$ (since $y<t$ and $z<t$); this suggests breaking into four cases: $0\le y+3z<t$, $t\le y+3z<2t$, $2t\le y+3z<3t$, and $3t\le y+3z<4t$. These cases can be analyzed one by one; for example, if $2t\le y+3z<3t$, then
\[t(x+y+z+2)+(y+3z-2t)=at+b\]
is a base $t$ expression, meaning $a=x+y+z+2$ and $b=y+3z-2t$. Combining all of this together, the set $S_t\cap[0,t^2)$ is in bijection with the set of $(a,b)$ (under the map $(a,b)\mapsto at+b$)  such that
\begin{align*}
\exists x,y,z\in\Z,\ &(0\le a,b,x,y,z< t)\\
&\wedge\ \bigvee_{i=0}^{3}\Big((a=x+y+z+i)\ \wedge \ (b=y+3z-it)\Big).
\end{align*}
This is getting too ugly to want to continue analyzing, but we could. No variables in this expression are multiplied by $t$! (Remember that $i$ is just one of the constants $0,1,2,3$.) The next step would be to eliminate the quantifiers, and we will find $F(t,t+1,t+3)$ and $g(t,t+1,t+3)$ as EQP functions. See the teaser in Section 1 for the final answer.
\end{example}

In addition to the Frobenius number and the number of gaps, many other invariants of parametric semigroups have recently been studied; see \cite{garcia2020,garcia2019b,garcia2019factorization}, for example. The connection of these \emph{parametric semigroups} to the methods we are describing here was explored in \cite{bogart2019periodic}.

\section{Statement of Theorems.}
\label{sec:thms}

In this section, we state everything precisely. Recall the definitions of quasi-polynomial (QP) and piecewise quasi-polynomial (PQP) functions of several parameters (Definitions \ref{def:qp} and \ref{def:pqp}). Recall that, with a single parameter $t$, PQP is more conveniently conceived of as \emph{eventually quasi-polynomial} (EQP): that is, agreeing with a quasi-polynomial for sufficiently large values of $t$. 

\begin{definition}
We say that a subset $T\subseteq\Z^k$ is \emph{piecewise periodic} if
\[\chi_T(\vec t)=\begin{cases}1&\text{if $\vec t\in T$,}\\ 0&\text{else}\end{cases}\]
is PQP, that is, if $T$ is a finite union of intersections of lattice cosets and polyhedra. With $k=1$, this simply means that the set $T$ is eventually periodic.
\end{definition}

In general, let $S_{\vec t}$, for ${\vec t}\in\Z^k$, be a family of subsets of $\Z^d$ and consider the following properties (of a quasi-polynomial nature) that $S_{\vec t}$ might have.

\begin{description}
\item[Property 1.] The set of ${\vec t}$ such that $S_{\vec t}\ne \emptyset$ is piecewise periodic.
\end{description}

\begin{description}
\item[Property 2.] There exists a PQP $g:\Z^k\rightarrow\Z$ such that, if $S_{\vec t}$ has finite cardinality, then $g({\vec t})=\abs{S_{\vec t}}$. The set of ${\vec t}$ such that $S_{\vec t}$ has finite cardinality is piecewise periodic.\end{description}

\begin{description}
\item[Property 3.] There exists a function $\x:\Z^k\rightarrow\Z^d$, whose coordinate functions are PQPs, such that, if $S_{\vec t}$ is non\-emp\-ty, then $\x({\vec t})\in S_{\vec t}$. The set of ${\vec t}$ such that $S_{\vec t}$ is nonempty is piecewise periodic.
 \end{description}

\begin{description}
\item[Property  3a.] Given $\vec c\in \Z^d\setminus\{0\}$, there exists a function $\vec x:\Z^k\rightarrow\Z^d$, whose coordinate functions are PQPs, such that if $\max_{\y\in S_{\vec t}} \vec c\cdot \y$ exists, then it is attained at $\x({\vec t})\in S_{\vec t}$. The set of  ${\vec t}$ such that the maximum exists is piecewise periodic.
\end{description}

\begin{description}
\item[Property  3b.] Fix $n\in\Z_{>0}$. There exist functions $\vec x_1,\ldots,\vec x_n:\Z^k\rightarrow\Z^d$, whose coordinate functions are PQPs, such that if $\abs{S_{\vec t}}\ge n$, then $\x_1(t),\ldots,\x_n({\vec t})$ are distinct elements of $S_{\vec t}$. The set of ${\vec t}$ such that $\abs{S_{\vec t}}\ge n$ is piecewise periodic.
\end{description}

%
%
%
Property 2 is designed for counting problems, while Properties 3, 3a, and 3b apply to various types of specification problems. Property 1 is a more fundamental property about nonemptiness, important in its own right, but also needed for the other definitions.  

Now we define some families, $S_{\vec t}$, that we might hope have these properties. These families are crafted using \emph{quantification over integers} ($\exists x\in\Z$ or $\forall x\in\Z$, where $x$ is one of the nonparameter variables, $\vec x$) and \emph{boolean combinations} (stringing together simpler formulas using the logical operations and, or, not). First, consider the sets with multiple parameters that we had defined up through Section \ref{sec:quant}:

\begin{definition}
A \emph{Presburger family} is a collection $\{S_{\vec t} : {\vec t} \in \Z^k\}$ of subsets of $\Z^d$ that can be defined by quantifying a boolean combination of formulas of the form $\vec a\cdot \x \le b+\vec c\cdot\vec t,$ where $b,\vec a,\vec c$ are constants with integer coordinates.
\end{definition}

Next, consider the family of sets where we restricted to one parameter variable, but we allowed nonlinearity in that parameter, as in Section \ref{sec:nonlinearity}:
\begin{definition}
A \emph{1-parametric Presburger family} is a collection $\{S_t : t \in \Z\}$ of subsets of $\Z^d$ that can be defined by quantifying a boolean combination of formulas of the form $\vec a(t)\cdot \x \le b(t),$ where  $b(t),\vec a(t)$ have coordinates in $\Z[t]$.
\end{definition}

Our two main theorems say that these families have everything we could possibly want!

\begin{theorem}
A Presburger family $S_{\vec t}$ has all of the Properties 1, 2, 3, 3a, and 3b.
\end{theorem}
This is proved in \cite{woods2015presburger}; Properties 3/3a/3b do not appear there, but they follow immediately from the techniques presented in \cite{Woods_unreasonable}.

\begin{theorem}
A 1-parametric Presburger family $S_{t}$ has all of the Properties 1, 2, 3, 3a, and 3b.
\end{theorem}
This is proved in \cite{BGW2017}. Both Presburger families and 1-parametric Presburger families also have certain generating function properties that we will not discuss here; see \cite{BGW2017,Woods_unreasonable,woods2015presburger}.

\section{Going too far.}
Where to next? It's not so clear. Here are two options for generalizing that lead to problems: allow a second (nonlinear) parameter or allow nonlinearity in the nonparameter variables.

\begin{figure}
\begin{center}
\includegraphics[width=.5\textwidth]{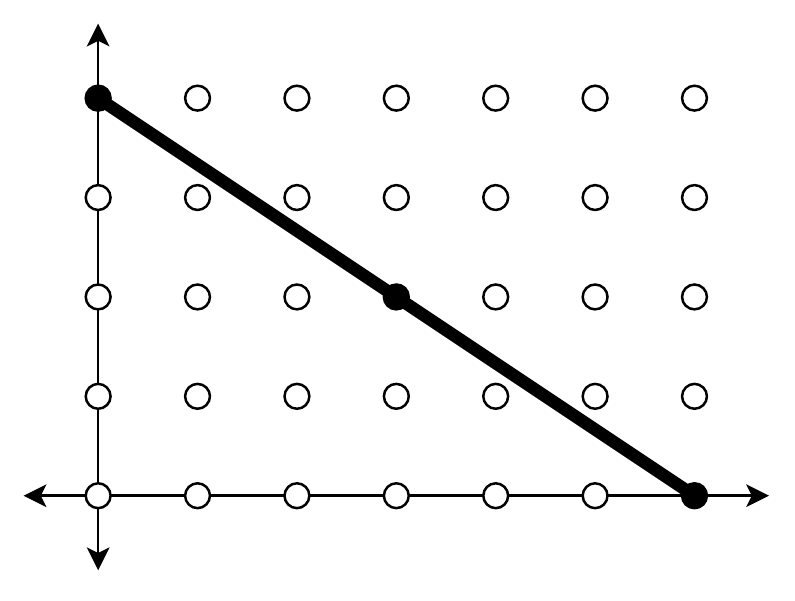}
\caption{Taking $s=4$ and $t=6$ in Example \ref{ex:gcdpic}, we examine the line segment of $x,y\ge 0$ such that $4x+6y=4\cdot 6$, which has endpoints $(0,4)$ and $(6,0)$. Therefore $f(4,6)=3=\gcd(4,6)+1$.}
\label{fig10}
\end{center}
\end{figure}

\begin{example} \cite{Woods_unreasonable}
\label{ex:gcdpic}
Let $f(s,t)$ be the number of $(x,y)\in\Z_{\ge 0}^2$ such that $sx+ty=st$. This has two nonlinear parameters, $s$ and $t$,  and it defines the interval in $\Z^2$ with endpoints $(0,s)$ and $(t,0)$; see Figure \ref{fig10}. So $f(s,t)=\gcd(s,t)+1$, which is a very nonpolynomial function. 
\end{example}

\begin{example} \cite{Woods_unreasonable}
Let $f(t)$ be the number of $x\in\Z_{\ge 0}$ such that $\exists y\in\Z_{\ge 0},\ xy=t$. This is nonlinear in $x$ and $y$, and $f(t)$ is the number of divisors of $t$, which is a very nonpolynomial function.
\end{example}

We should step back and be thankful that nonlinearity in a parameter variable ended up so nicely (we will discuss a little bit of why we think it ended up so nicely in Section \ref{sec:conc}). Allowing nonlinearity is extremely dangerous, so be careful! Allowing such expressiveness rapidly approaches the domain of G\"odel's first incompleteness theorem:

\begin{example}
Let $p(x_0,\ldots,x_d)$ be a polynomial, and define $f_p(t)$ to be the number of $(x_1,\ldots, x_d)\in\Z^d$ such that $p(t,x_1,\ldots,x_d)=0$. A consequence of the DPRM theorem (see \cite{Davis73}), which solves Hilbert's 10th problem in the negative, is that there exists a specific polynomial $p$ such that $f_p(t)$ is not computable (and indeed, that figuring out whether $f_p(t)=0$ is undecidable.) 
\end{example}

\section{Not Far Enough?}

\begin{figure}
\begin{center}
\includegraphics[width=.4\textwidth]{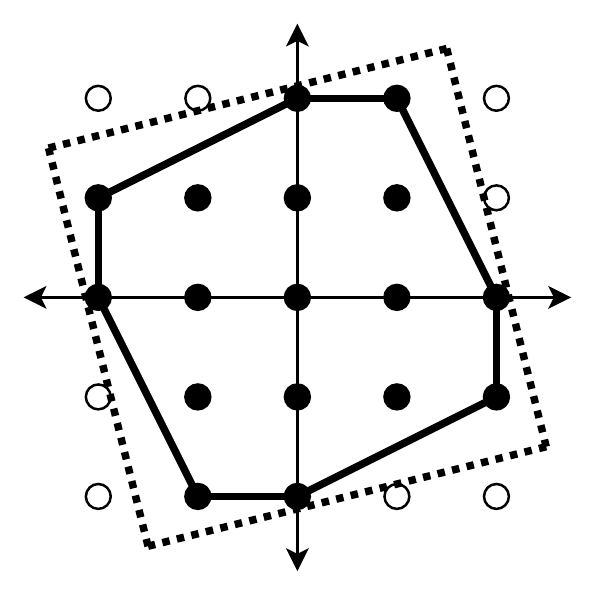}
\caption{The integer hull (solid line) of $P_5$ (dotted line), in Example \ref{ex:twist} (cf., Figure \ref{fig9}). When $t$ is odd, the integer hull of $P_t$ is an octagon with vertices $\left(0,\pm\frac{t-1}{2}\right),\left(\pm\frac{t-3}{2},\pm\frac{t-1}{2}\right),\left(\pm\frac{t-1}{2},0\right),\left(\pm\frac{t-1}{2},\mp\frac{t-3}{2}\right).$ When $t$ is even, the vertices of $P_t$ are integers, so the vertices of the integer hull are simply the vertices of $P_t$:
 $\left(\pm\frac{t-2}{2},\pm\frac{t}{2}\right)$ and $\left(\pm\frac{t}{2},\mp\frac{t-2}{2}\right)$; see Example \ref{ex:ih}.}.
\label{fig11}
\end{center}
\end{figure}

But we shouldn't give up so easily. One might look for other natural problems that, though they do not seem to fit into the domain of parametric Presburger arithmetic, do still yield EQP results:
\begin{example}
\label{ex:ih}
For $t\in\Z_{\ge 0}$, let $P_t\subseteq\Z^d$ be a polytope defined by linear equations of the form $\vec a(t)\cdot\vec x\le b(t)$, where $\vec a:\Z\rightarrow \Z^d$ and $b:\Z\rightarrow\Z$ have polynomial coordinates.  It has been shown that there is an EQP function $\vec v:\Z_{\ge 0}\rightarrow \Z^d$ such that $\vec v(t)$ is a vertex of the convex hull of $P_t\cap\Z^d$ \cite{CalegariWalker, Shen2018}; see, for example, Figure \ref{fig11}. A natural way to specify such a $\vec v$ would be 
\[ (\vec v \in P_t\cap \Z^d) \wedge \exists \vec z\in\Z^d,\ \forall \vec x\in P_t\cap \Z^d, \ \left( \vec x = \vec v\right)\ \vee\ \left(\vec z\cdot\vec x < \vec z\cdot\vec v \right). \]
  Note that writing $\vec z\cdot\vec x$ and $\vec z\cdot\vec v$ puts us outside of the domains of Presburger arithmetic and 1-parametric Presburger arithmetic, and yet we still get an EQP specification function.
\end{example}

\begin{example}
For $t\in\Z_{\ge 0}$, let $\Lambda_t\subseteq\Z^d$ be the lattice generated by $\vec b_1(t),\ldots,\vec b_d(t)$, where the $\vec b_i$ have polynomial coordinates. Then there is an EQP function $\vec v:\Z_{\ge 0}\rightarrow \Z^d$ such that $\vec v(t)$ is the shortest (in the Euclidean norm) nonzero vector in $\Lambda_t$ \cite{bogart2019parametric}. Comparing Euclidean norms of vectors seems to require nonlinearity, since $||\vec v||^2=\sum_i v_i^2$. So again we seem to have left the domains of Presburger arithmetic and 1-parametric Presburger arithmetic, and yet we still get an EQP function.
\end{example}

In both of these examples, additional tools (specific to each problem) are required to prove EQP properties. 
\begin{question} Is there a common generalization: some larger class of problems whose answer can be shown to be EQP?
\end{question}

Another way to generalize might be to allow two nonlinear parameters, accept that we won't get EQP functions, but still hope that our functions are ``nice'' somehow. For example, the lcm function pops up in some simple examples:

\begin{example}
Let's further define the Frobenius problem for generators that are not relatively prime. Let $F(a_1,\ldots,a_d)$ be the largest element of $a_1\Z+\cdots+a_d\Z$ that is not in the semigroup generated by $a_1,\ldots,a_d$. Then
\[F(a_1,a_2)=\lcm(a_1,a_2)-a_1-a_2.\]
\end{example}

What a nice function!

\begin{question}Can we generalize to other 2-parametric Presburger families?
\end{question}

In full generality, the answer is an emphatic ``No!''; not all 2-parametric families have cardinalities that are nice functions. The proof of this acknowledges that deciding what functions are ``nice'' is hard, but asks us to believe that nice functions can at least be computed in polynomial time (in a theoretical computer science sense; polynomials, QPs, and lcms, for example, can be evaluated in polynomial time). Assuming  P $\ne$ NP, \emph{there exists a specific 2-parametric family $S_{s,t}$ whose cardinality cannot be computed in polynomial time in the size of the input $s$ and $t$} \cite{BGNW}.

On the other hand, if we are only allowed equalities rather than inequalities, then it turns out that functions related to lcm, gcd, and the extended Euclidean algorithm do suffice to describe all counting functions, for any number of nonlinear parameters \cite{BGNW} (as usual, the nonparameter variables must be linear, or we would run into G\"odel-like problems quickly).

\section{Why EQP?}
\label{sec:conc}

To conclude, we would like to philosophize about why piecewise/eventual quasi-polynomials show up in so many counting problems related to discrete geometry.

First of all, why not simply polynomials? Well, even a simple formula like $0\le 2x\le t$ immediately requires some periodicity, with functions like $\floor{t/2}$, so we would expect quasi-polynomials to appear fairly often.

So why not simply QPs? Well, with one parameter, simple formulas like $0\le (3t+2)x\le t^2$ will require us to calculate floor functions, like $\floor{t^2/(3t+2)}$, that are only quasi-polynomial for sufficiently large $t$. Or we might have to decide which of two polynomials, like $t^2+10$ and $2t^2-5$, is larger, which will depend on whether $t$ is sufficiently large. So we would expect EQPs to appear fairly often. With more than one parameter, we immediately see examples like the trapezoid/triangle phenomenon in Example \ref{ex:trap}; the combinatorial structure changes as the parameters change, and so we need piecewise quasi-polynomials.

The bigger question is why EQP functions suffice for so many counting problems. This seems to follow directly from the fact that the floor of a rational function, $\floor{f(t)/g(t)}$, is EQP. The floor function is the basis of many operations that we might need to undertake, such as the gcd operation. Algebraically, any finitely generated ideal in the ring of EQPs is principal \cite{CLS} (such a ring is called a \emph{Bezout domain}). While the ring of integer-valued polynomials and the ring of QPs are not closed under some very natural operations, the ring of EQPs is.

This also helps us see why problems with two parameters can get rapidly hard.

\begin{question}What kind of function is $\floor{\frac{t^2+st+3s+7}{t-s}}$? Can it be simplified in any reasonable way, in the sense that $\floor{f(t)/g(t)}$ can be simplified to an EQP? We invite the reader to explore!
\end{question}


\begin{acknowledgment}{Acknowledgments.} We thank John Goodrick, Danny Nguyen, and Chris O'Neill for useful discussions as we were writing this paper. We also appreciate the reviewers' careful reading and helpful comments. 
\end{acknowledgment}

\begin{biog}
\item[Tristram Bogart] received his Ph.D. from the University of Washington in 2007. After post docs at Queen's University and MSRI / San Francisco State University, he moved to Los Andes University in Bogot\'a, Colombia, where he is now an Associate Professor of Mathematics. 
\begin{affil}
Departamento de Matem\'aticas, Universidad de los Andes, Bogot\'a, Colombia\\
tc.bogart22@uniandes.edu.co
\end{affil}

\item[Kevin Woods] received his Ph.D. from the University of Michigan in 2004. After a post doc at the University of California, Berkeley, he moved to Oberlin College in 2006, and he is now a Professor of Mathematics there. He is happiest when he is with his dogs or is taking long walks, but he is sad that his dogs cannot keep up on these long walks.
\begin{affil}
Department of Mathematics, Oberlin College, Oberlin, OH 44074\\
Kevin.Woods@oberlin.edu
\end{affil}
\end{biog}
\vfill\eject

\end{document}